\documentclass[fleqn,10pt]{wlscirep}
\usepackage[utf8]{inputenc}
\usepackage{times}
\usepackage{eurosym}
\usepackage{titlesec}
\usepackage{amsmath}
\usepackage{amscd, fullpage,amssymb, color,amsthm, todonotes}
\usepackage[all]{xy}
\usepackage{float}
\usepackage{blkarray}
\usepackage[T1]{fontenc}
\usepackage{xcolor}
\usepackage{comment}
\titleformat{\paragraph}
{\normalfont\normalsize\bfseries}{\theparagraph}{1em}{}
\titlespacing*{\paragraph}
{0pt}{3.25ex plus 1ex minus .2ex}{1.5ex plus .2ex}

\theoremstyle{remark}

\theoremstyle{definition}
\newtheorem{defn}{Definition}

\definecolor{brilliantrose}{rgb}{1.0, 0.33, 0.64}
\definecolor{myviolet}{rgb}{0.21, 0.0, 0.85}
\definecolor{amethyst}{rgb}{0.6, 0.4, 0.8}
\definecolor{carrotorange}{rgb}{0.93, 0.57, 0.13}

\title{A cohomology-based Gromov-Hausdorff metric approach for quantifying molecular similarity}

\author[1*]{JunJie Wee}
\author[2$\dag$]{Xue Gong}
\author[3$\ddagger$]{Wilderich Tuschmann}
\author[2$\mathsection$]{Kelin Xia}
\affil[1]{Michigan State University, Department of Mathematics, East Lansing, MI 48824, USA}
\affil[2]{Nanyang Technological University, Division of Mathematical Sciences, School of Physical and Mathematical Sciences, 637371, Singapore}
\affil[3]{Karlsruher Institut für Technologie, Institut für Algebra und Geometrie, Arbeitsgruppe Differentialgeometrie, Karlsruhe, 76131, Germany}
\affil[*]{weejunji@msu.edu} 
\affil[$\dag$]{xue.gong@ntu.edu.sg}
\affil[$\ddagger$]{wilderich.tuschmann@kit.edu}
\affil[$\mathsection$]{xiakelin@ntu.edu.sg}

\begin{abstract}

We introduce, a cohomology-based Gromov-Hausdorff ultrametric method to analyze 1-dimensional and higher-dimensional (co)homology groups, focusing on loops, voids, and higher-dimensional cavity structures in simplicial complexes, to address typical clustering questions arising in molecular data analysis. The Gromov-Hausdorff distance quantifies the dissimilarity between two metric spaces. In this framework, molecules are represented as simplicial complexes, and their cohomology vector spaces are computed to capture intrinsic topological invariants encoding loop and cavity structures. These vector spaces are equipped with a suitable distance measure, enabling the computation of the Gromov-Hausdorff ultrametric to evaluate structural dissimilarities.  We demonstrate the methodology using organic-inorganic halide perovskite (OIHP) structures. The results highlight the effectiveness of this approach in clustering various molecular structures. By incorporating geometric information, our method provides deeper insights compared to traditional persistent homology techniques.

 \end{abstract}

\begin{document}

\flushbottom
\maketitle

%\keywords{Keyword1, Keyword2, Keyword3}
%  Click the title above to edit the author's information and abstract
%
\thispagestyle{empty}

\section*{Introduction}
In the last decade, concepts, methods, and techniques from topological and geometric data analysis \cite{carlsson2009topology, wasserman2018topological, le2004geometric, kirby2000geometric} have found profound applications, becoming indispensable tools in many areas of applied mathematics, such as image processing \cite{klette2004digital},  pattern recognition \cite{carlsson2014topological}, sensor network \cite{gao2012geometric}, robotics \cite{bhattacharya2010search}, computer graphics and vision \cite{skraba2010persistence}, cosmology \cite{sousbie2011persistent},  medicine \cite{skaf2022topological}, computational and molecular biology \cite{yao2009topological}. A central technique in Topological Data Analysis (TDA) is persistent homology, which provides a robust framework for capturing multi-scale topological features in data \cite{edelsbrunner2002topological, zomorodian2004computing}. However, while persistent homology is a powerful tool, it generally provides only limited information about the homology groups of the objects under study, often restricted to basic descriptors such as Betti numbers. In addition to topological tools, geometric theories and methods, such as discrete exterior calculus \cite{hirani2003discrete}, Laplace operators \cite{belkin2003laplacian}, discrete Ricci curvatures \cite{ollivier2009ricci, forman2003bochner}, optimal transport \cite{villani2009optimal}, parametric surfaces \cite{farin2014curves}, and geometric flows \cite{jin2008discrete}, have been widely applied in data analysis.  These geometric approaches have often been shown to provide deeper insights than those achieved through purely topological methods. 

A promising application of advanced geometric and topological techniques is molecular similarity analysis \cite{bender2004molecular}, which is a fundamental concept in cheminformatics. Molecular similarity measures the resemblance between molecules based on their structure, properties, or function, enabling applications in drug design \cite{vilar2012drug}, virtual screening \cite{eckert2007molecular}, and reaction prediction \cite{arteca1988shape}. By assessing similarity, we can predict properties such as activity and toxicity, thereby accelerating the identification of promising compounds \cite{petrone2012rethinking, basak1998comparative}.

Existing methods for evaluating molecular similarity often rely on structural fingerprints, property-based metrics, or topological and 3D shape representations. For instance, structural fingerprints, such as Extended Connectivity Fingerprint (ECFP) \cite{rogers2010extended} and MACCS keys \cite{polton1982installation}, encode molecular features as binary sequences, which are typically compared using Tanimoto similarity \cite{bajusz2015tanimoto}. Property-based similarity involves comparing physicochemical attributes, such as molecular weight and polarity, using Euclidean or cosine distances \cite{todeschini2008handbook}. Topological similarity, often evaluated through graph neural networks (GNNs) \cite{gilmer2017neural}, models molecular connectivity, while 3D similarity captures spatial geometry through metrics such as RMSD or Shape Tanimoto scores \cite{grant1996fast, good1993rapid}.

Recent advances in applied and computational topology and geometry offer powerful methods to assess molecular similarity. Persistent homology and its various refinements and extensions have proven invaluable for capturing molecular topological features in diverse chemical and biological applications \cite{townsend2020representation, xia2014persistent}. However, while effective in tracking the number of connected components, loops, and voids, persistent homology lacks the geometric information necessary to distinguish or identify individual features that arise during filtration. Previous studies have leveraged the connection between persistent homology and dendrograms \cite{carlsson2010characterization, carlsson2008persistent} to distinguish structures based on 0-dimensional homology groups, which correspond to connected components. These works employ geometric methods, such as the Gromov-Hausdorff distance, to measure the distance between two brain networks \cite{lee2011computing, lee2012persistent, chung2017topological}. However, there is still a need for methodologies that enable a more detailed and localized characterization of loops, voids, and higher-dimensional (co)homologies.

We apply a cohomology-based Gromov-Hausdorff ultrametric method to study 1-dimensional and higher-dimensional (co)homology groups associated with loops, voids, and higher-dimensional cavity structures in simplicial complexes. The Gromov-Hausdorff distance ($d_{GH}$) \cite{edwards1975structure, gromov1981groups} quantifies the dissimilarity between two metric spaces and possesses many desirable mathematical properties, making it widely used in differential geometry and applied algebraic topology for tasks such as shape matching, structural comparison, and characterization \cite{memoli2007use, bronstein2010gromov, chazal2009gromov}.  In the context of molecular structures $A$ and $B$, modeled as simplicial complexes $K_A$ and $K_B$, respectively, their cohomology groups give rise to the metric spaces $H^p(K_A)$ and $H^p(K_B)$. In theory, we can then compare their Gromov-Hausdorff distances. Although many lower bounds can be achieved in polynomial time \cite{chazal2009gromov, memoli2012some}, computing $d_{GH}$ between arbitrary finite metric spaces is known to be NP-hard, as it leads to quadratic assignment problems \cite{memoli2012some}. To address this challenge, we introduce a crucial second step: transforming $H^p(K_A)$ and $H^p(K_B)$ into ultrametric spaces, thus obtaining their respective dendrogram representations \cite{carlsson2010characterization, memoli2021gromov, lee2012persistent}. This enables us to compute a more tractable ultrametric variation of the Gromov-Hausdorff distance, referred to as the Gromov-Hausdorff ultrametric ($u_{GH}$)\cite{zarichnyi2005gromov, memoli2021gromov}, which can be used to quantify molecular similarity.

The main contributions of this work are as follows: 
\begin{enumerate} 
\item We present a workflow to utilize a Gromov-Hausdorff ultrametric approach based on 1-dimensional and higher-dimensional cohomology for assessing structural similarity between molecules, effectively capturing local topological features such as loops and voids that influence molecular properties.

\item We demonstrate the effectiveness of this workflow in clustering tasks using a dataset of organic-inorganic halide perovskites (OIHP). 
\end{enumerate}

The remainder of this paper is structured as follows: The next section introduces foundational concepts in algebraic topology, including simplicial complexes and the Hodge Laplacian, which are essential to our methodology. We then present our method for measuring structural similarity, followed by numerical experiments to validate our approach using OIHP material datasets. Finally, we discuss potential future research directions.

%%%%%%%%%%%%%%%%%%%%%%%%%%%%%%%%%%%%%%%%%%%%%%%%%%%%%%%%%%%%%%%%%%%%%%%%%%%%%%
\section*{Background}
%----\textcolor{orange}{no metric yet, purely combinatorial}}
\subsection*{Simplicial complexes}

In many applications, molecular data from biology, chemistry, and material science can be represented as graphs \cite{bajorath2004chemoinformatics, lo2018machine}. In this framework, the vertices correspond to atoms, while the edges represent affinities between pairs of atoms. Simplicial complexes generalize graphs by capturing higher-order interactions through simplices, such as triangles and tetrahedra. In this and the next section, we aim to provide a brief introduction to the theoretical background necessary for performing computations on simplicial complexes. For a more comprehensive review, we direct readers to \cite{eckmann1944harmonische,lim2015hodge,shukla2020spectral,torres2020simplicial}.

If we let $V$ be a finite set of vertices, then a $p$-simplex $\sigma^p$ is a subset of $V$ with $p+1$ vertices. For example, a simplex of dimensions 0, 1, 2, and 3 can be viewed as a point, an edge, a triangle, and a tetrahedron, respectively. More precisely, a $p$-simplex $\sigma^p = \{v_0, v_1, v_2, \cdots, v_p\}$ is defined as a convex hull formed by its $p+1$ affinely independent points $v_0, v_1, v_2, \cdots, v_p$ as follows:
\[
\sigma^p = \bigg\{\lambda_0v_0+\lambda_1v_1+\cdots+\lambda_pv_p \bigg|\sum_{i=0}^p \lambda_i = 0;\forall i, 0\leq \lambda_i \leq 1\bigg\}.
\]
A $m$-face of a $p$-simplex $\sigma^p$ is defined as a convex hull formed by $m+1$ vertices of $\sigma^p$, where $m < p$. If $\sigma^{m}$ is a face of $\sigma^p$, denoted by $\sigma^{m} \subset \sigma^p$, then $\sigma^p$ is also referred to as a \textit{coface} of $\sigma^{m}$. The \textit{upper degree} of a $p$-simplex $\sigma^p$, denoted by $deg(\sigma^p)$, is the number of $(p+1)$-simplices for which $\sigma^p$ is a face. Two simplices $\sigma_1^p$ and $\sigma_2^p$ are \textit{upper adjacent} and denoted by $\sigma_1^p \frown \sigma_2^p$ if they have a common coface. They are \textit{lower adjacent} and denoted by $\sigma_1^p \smile \sigma_2^p$ if they share a common face. To facilitate computation, each simplex is assigned an orientation associated with the ordering of its vertices. If two $p$-simplices $\sigma_1^p$ and $\sigma_2^p$ are oriented similarly, then we write $\sigma_1^p \sim \sigma_2^p$.  Conversely, we write $\sigma_1^p \not\sim \sigma_2^p$ if they have opposite orientations.

A $p$-dimensional simplicial complex $K$ contains up to $p$-dimensional simplices and satisfies two conditions. First, any face of a simplex from $K$ is also in $K$. Second, the intersection of any two simplices in $K$ is either empty or a shared face. There are many common types of simplicial complexes such as Vietoris-Rips complex, $\check{C}ech$ complex, Alpha complex, and clique complex. 

\subsection*{The Combinatorial Hodge Laplacian}

Computationally, combinatorial Hodge Laplacian or discrete Hodge Laplacian has been proposed \cite{eckmann1944harmonische,muhammad2006control,horak2013spectra,barbarossa2020topological,mukherjee2016random,parzanchevski2017simplicial,shukla2020spectral,torres2020simplicial}. Essentially, this discrete version can be viewed as part of exterior calculus and discrete differential geometry. The concept of the combinatorial Hodge Laplacian is a part of Hodge theory, which has recently been applied to biomolecular data analysis \cite{wei2022hodge}. Hodge Laplacian matrices of different dimensions can be constructed on a simplicial complex. A $p$-th dimensional Hodge Laplacian matrix characterizes topological connections between $p$-th simplices within a simplicial complex \cite{eckmann1944harmonische, lim2015hodge, doi:10.1137/18M1201019}. Note that the graph Laplacian characterizes relations between vertices (0-simplices).

The $p$-th chain group $C_p(K)$ of a simplicial complex $K$ over some field $\mathbb{F}$ is a vector space over $\mathbb{F}$ whose basis is the set of $p$-simplices of the simplicial complex $K$. Elements of $C_p(K)$ are called $p$-chains. The \textit{dual} of $C_p(K),$ denoted by $C^p(K),$ is the set of all linear functionals on $C_p(K)$:
$$C^p(K)=\big\{\phi: C_p(K) \to \mathbb{F}\, : \, \text{$\phi$ is linear} \big\}.$$
$C^p(K)$ is called the \textit{$p$-th cochain group} and its elements are called \textit{$p$-cochains}. Boundary operators are defined on both the chain and cochain groups. The \textit{boundary map} $\partial_p$ is a linear transformation which acts on a $p$-simplex $\sigma^p=[u_0,u_1,\ldots,u_p]$ as follows
$$\partial_p([u_0,u_1,\ldots,u_p])=\sum_{i=0}^p(-1)^i[u_0,\ldots,u_{i-1},u_{i+1},\ldots,u_p].$$
The \textit{coboundary map} $\delta_p:C^p(K) \to C^{p+1}(K)$ is a linear transformation defined as follows: for a linear functional $\phi \in C^p(K)$ and a $p+1$-simplex $\sigma^{p+1}=[u_0,u_1,\ldots,u_{p+1}]$,
$$\delta_p(\phi)(\sigma^{p+1})=\sum_{i=0}^{p+1}(-1)^i\phi([u_0,\ldots,u_{i-1},u_{i+1},\ldots,u_{p+1}]).$$
The boundary map gives rise to a \textit{chain complex}, which is a sequence of chain groups connected by boundary maps as follows:
$$ 0 \to C_{n}(K)\to\cdots\xrightarrow{\partial_{p+1}}C_p(K)\xrightarrow{\partial_p} C_{p-1}(K)\cdots\xrightarrow{\partial_2}C_1(K)\xrightarrow{\partial_1}C_0(K)\rightarrow 0.$$
Similar to the boundary map giving rise to the chain complex, the coboundary operator gives rise to a \textit{cochain complex}:
$$ 0 \leftarrow C^{n}(K)\leftarrow\cdots\xleftarrow{\delta_{p}}C^p(K)\xleftarrow{\delta_{p-1}} C^{p-1}(K)\cdots\xleftarrow{\delta_1}C^1(K)\xleftarrow{\delta_0}C^0(K)\leftarrow 0.$$
Since $C_p(K)$ and $C^p(K)$ are finite-dimensional, there exists unique matrix representations for $\partial_p$ and $\delta_p$. We have some useful relations regarding matrix representations of $\partial_p$ and $\delta_p$ ($A^T$ represents the transpose of a matrix $A$):
\begin{itemize}
	\item For all $p \geq 0$, $\partial_{p+1}^T=\delta_p$,
	\item $\partial_p^T=\partial_p^*$,
	\item $\delta_p^T=\delta_p^*.$
\end{itemize}
Here, $\delta_p^*:C^{p+1}(K)\to C^p(K)$ is the \textit{adjoint/transpose map} of $\delta_p$ where
$$\langle \delta_p(f), g \rangle=\langle f, \delta_p^*(g) \rangle,$$
for every $f \in C^p(K)$, $g \in C^{p+1}(K)$ and a suitable inner product $\langle \cdot , \cdot \rangle$ for $C^p(K)$ and $C^{p+1}(K).$ The adjoint of the boundary operator $\partial_p$, $\partial_p^*$ is also defined analogously. Using the cochain group $C^p(K)$, we can define the $p$-th cocycle group $Z^p$ and $p$-th coboundary group $B^p$ as follows:
\begin{equation*}
Z^p = \text{ker}(\delta_p) = \{c \in C^p | \delta_p(c) = 0\},
\end{equation*}
\begin{equation*}
B^p = \text{im}(\delta_{p-1}) = \{c\in C^p | \exists d \in C^{p-1} : c = \delta_{p-1}(d) \}.
\end{equation*}Then we have the $p$-th cohomology group $H^p = Z^p/B^p$.

The $p$-dimensional \textit{combinatorial Hodge Laplacian} is the linear operator $\Delta_p:C^p(K) \to C^p(K)$ is defined as follows:
$$\Delta_p=
\begin{cases}
\delta_p^*\circ\delta_p+\delta_{p-1}\circ\delta_{p-1}^* & \text{if } p \geq 1, \\
\delta_p^*\circ\delta_p & \text{if } p=0.
\end{cases}
$$
The case where $p=0$ gives rise to the expression of the well-known graph Laplacian. Alternatively, the combinatorial Hodge Laplacian matrix can be expressed using the boundary matrix. The boundary operator $\partial_p$ has a unique matrix representation. Given a simplicial complex $K$, the \textit{$p$-th boundary matrix} $\mathbf{B}_p$ is defined as,
$$(\mathbf{B}_p)_{ij}=
\begin{cases}
1 &\text{if } \sigma_i^{p-1} \subset \sigma_j^{p} ~\text{and }~ \sigma_i^{p-1} \sim \sigma_j^{p},  \\
-1 &\text{if } \sigma_i^{p-1} \subset \sigma_j^{p} ~\text{and }~ \sigma_i^{p-1} \not\sim \sigma_j^{p},  \\
0 &\text{if } \sigma_i^{p-1} \not\subset \sigma_j^{p},
\end{cases}
$$
where $\sigma_i^{p-1}$ is the $i$-th {$(p-1)$-simplex} and $\sigma_j^{p}$ is the $j$-th $p$-simplex.

Given that the highest order of the simplicial complex $K$ is $n$, the $p$-th {Hodge} Laplacian (or combinatorial Laplacian) matrix $\mathbf{L}_p$ of $K$ is
$$\mathbf{L}_p=
\begin{cases}
\mathbf{B}_n^T\mathbf{B}_n & \text{if } p=n, \\
\mathbf{B}_p^T\mathbf{B}_p+\mathbf{B}_{p+1}\mathbf{B}_{p+1}^T & \text{if } 1 \leq p < n, \\
\mathbf{B}_1\mathbf{B}_1^T & \text{if } p=0.
\end{cases}
$$

Another way to define the Hodge Laplacian is through simplex relations. When $p=0$,
$$(\mathbf{L}_0)_{ij}=
\begin{cases}
\text{deg}(\sigma_i^{0}) &\text{if } i=j, \\
-1 & \text{if } i \neq j~\text{and } \sigma_i^{0} \frown \sigma_j^{0}, \\
0 & \text{if } i \neq j~\text{and } \sigma_i^{0} \not \frown \sigma_j^{0},
\end{cases}
$$
and $\mathbf{L}_0$ is exactly the same as the graph Laplacian matrix. When $p>0$,
$$(\mathbf{L}_p)_{ij}=
\begin{cases}
\text{deg}(\sigma_i^{p})+ p+1 \quad &\text{if } i=j, \\
1 \qquad &\text{if } i \neq j, \, \sigma_i^{p} \not \frown\sigma_j^{p}, \sigma_i^{p} \smile \sigma_j^{p}~\text{and } \sigma_i^k \sim  \sigma_j^p,   \\
-1 \quad &\text{if } i \neq j, \, \sigma_i^{p} \not \frown \sigma_j^{p}, \sigma_i^{p} \smile \sigma_j^{p}~\text{and } \sigma_i^k \not \sim  \sigma_j^p, \\
0 \qquad &\text{if } i \neq j, \, \sigma_i^{p} \frown \sigma_j^{p} ~\text{or }~ \sigma_i^{p} \not \smile \sigma_j^{p}.
\end{cases}
$$

Mathematically, the eigenvalues of Hodge Laplacian matrices are independent of the choice of the orientation \cite{horak2013spectra}. The eigenspectrum of Hodge Laplacian matrices reveals topological information within simplicial complexes. For instance, the multiplicity of zero eigenvalues of $\mathbf{L}_p$ corresponds to the Betti numbers $\beta_p$, reflecting the number of connected components ($\beta_0$), the number of cycles ($\beta_1$), the number of cavities ($\beta_2$), etc.

Furthermore, the eigenvectors corresponding to the zero eigenvalues of $\mathbf{L}_p$, also referred to as cohomology generators, can be used to represent the $p$-dimensional cohomologies. At the algebraic level, choosing a representative for an element in cohomology is inherently problematic due to its quotient structure. However, by leveraging the Hodge decomposition theorem with the inner product structure on the cochain complex, we obtain $\operatorname{ker}(\delta_p) = \operatorname{im}(\delta_{p-1}) \oplus \operatorname{ker}(\textbf{L}_p)$, where \( \operatorname{ker}(\textbf{L}_p) \) is orthogonal to \( \operatorname{im}(\delta_{p-1}) \). This identification induces the canonical isomorphism $H^p(K) \simeq (\operatorname{im}(\delta_{p-1}) \oplus \operatorname{ker}(\textbf{L}_p))/\operatorname{im}(\delta_{p-1}) \simeq \operatorname{ker}(\textbf{L}_p)$. In particular, every cohomology class in \( H^p(K) \) can be uniquely represented by an element in \( \operatorname{ker}(\textbf{L}_p) \) \cite{lim2015hodge, gurnari2023probing, basu2024harmonic}. For any cohomology generator $v$ of $\mathbf{L}_p$, we can arrange the entries of $v$ to match the order of simplices in $K$. Specifically, the $i$-th entry $v_i$ corresponds to the $p$-simplex $\sigma_i^p$ in $K$. This allows us to visualize $v$ on a simplicial complex. 

%v = \begin{blockarray}{cc}
%\begin{block}{(c)c}
%v_1 &  \sigma_{1}^p  \\
%v_2 &  \sigma_{2}^p \\
%\vdots &  \vdots \\
%v_{n_p-1} & \sigma_{n_p-1}^p  \\
%v_{n_p} &  \sigma_{n_p}^p \\
%\end{block}
%\end{blockarray},
%\]

%%%%%%%%%%%%%%%%%%%%%%%%%%%%%%%%%%%%%%%%%%%%%%%%%%%%%%%%%%%%%%%%%%%%%%%%%%%%%%%%%%%%%%%%%%%%%%%%

%%%%%%%%%%%%%%%%%%%%%%%%%%%%%%%%%%%%%%%%%%%%%%%%%%%%%%%%%%%%%%%%%%%%%%%%%%%%%%%%%
%\section*{\textcolor{orange}{After metrization, we get a metric space. Given two such metric spaces, how to compute the Gromov-%Hausdorff distance between them? It is usually difficult. Strategy: ultrametric them to get a dendrogram for each metric space, and then compute the ultra-Gromov-Hausdorff distance of the two dendrograms}}
%\end{comment}

\subsection*{The Gromov-Hausdorff Distance}

The Gromov–Hausdorff distance  measures the distance between two compact metric spaces \cite{edwards1975structure, gromov1981groups}.  It measures the smallest distance at which two compact metric spaces can be considered ``close'', which is particularly useful when the overall shape matching is more relevant than exact pointwise alignment. Given two metric spaces $X$ and $Y$, the Gromov–Hausdorff distance looks for all possible isometric embeddings of $X$ and $Y$ into a common metric space $Z$ and then calculates the Hausdorff distance \cite{hausdorff1914grundzuge} between these embeddings within $Z$. Formally, the Gromov-Hausdorff distance $d_{GH}$ \cite{edwards1975structure, gromov1981groups}  between two compact metric spaces $X$ and $Y$ is defined as 
\[
d_{GH}(X,Y) = \inf d_H^Z(\varphi(X),\psi(Y)),
\]
where the infimum is taken over all possible $Z\in \mathcal{M}$ and isometric embeddings $\varphi:X\hookrightarrow Z$ and $\psi:Y\hookrightarrow Z$, and $d^Z_H$ denotes Hausdorff distance in $Z$.

\section*{Method}

\subsection*{Cohomology Generator-based Distance Metrics}
For a given simplicial complex $K$, we denote the space consisting of all the $p$-dimensional cohomology generators arising from its Hodge Laplacian $\mathbf{L}_p$ by $H^p(K)$. Subsequently, we want to construct a metric space from  $H^p(K)$ by assigning distances between  cohomology generators in $H^p(K)$.  We will introduce three types of distance measures to construct the metric space: $L_1$ distance, cocycle distance, and Wasserstein distance.

Throughout this section, we assume that any two cohomology generators $v$ and $w$ from a Hodge Laplacian matrix $\mathbf{L}_p$ must have a consistent order of entries. In other words, we can write
\[
v = \begin{blockarray}{cc}
\begin{block}{(c)c}
v_1 &  \sigma_{1}^p  \\
v_2 &  \sigma_{2}^p \\
\vdots &  \vdots \\
v_{n_p-1} & \sigma_{n_p-1}^p  \\
v_{n_p} &  \sigma_{n_p}^p \\
\end{block}
\end{blockarray}
\text{ and }\quad
w = \begin{blockarray}{cc}
\begin{block}{(c)c}
w_1 &  \sigma_{1}^p  \\
w_2 &  \sigma_{2}^p \\
\vdots &  \vdots \\
w_{n_p-1} & \sigma_{n_p-1}^p  \\
w_{n_p} &  \sigma_{n_p}^p \\
\end{block}
\end{blockarray},
\]
where $v_i$ and $w_i$ are entries in $v$ and $w$, respectively, corresponding to the simplex $\sigma_i^p$ for all $1\leq i \leq n_p$. Furthermore, we assume that any cohomology generator in $H^p(K)$ is normalized. 

\subsubsection*{$L_1$ Distance}
First, we define the $L_1$ distance between two cohomology generators as follows. 
\begin{defn}[$L_1$ distance]\label{defn:pdist1}
Let $v$ and $w$ be two cohomology generators from $H^p(K)$ where $K$ is a simplicial complex. Then the $L_1$ distance between $v$ and $w$ is 
	\[
	\|v-w\|_1 = \sum_{i=1}^{n_p} |v_{i}-w_{i}|.
	\]
Note that in computation, the eigenvectors may have arbitrary signs. To eliminate the sign ambiguity introduced by the solver, we enforce that the first element of each cohomology generator is non-negative.  
\end{defn}

\subsubsection*{Cocycle Distance}
As the $L_1$ norm of a cohomology generator tends to be larger when the cocycle contains more edges, it serves as a rough indicator of the size of the cocycle. Therefore, it is natural to consider another type of distance that measures the absolute difference between the $L_1$ norms of two cohomology generators. We refer to this as the \textit{cocycle distance}.

\begin{figure}[ht]
	\includegraphics[width=\textwidth]{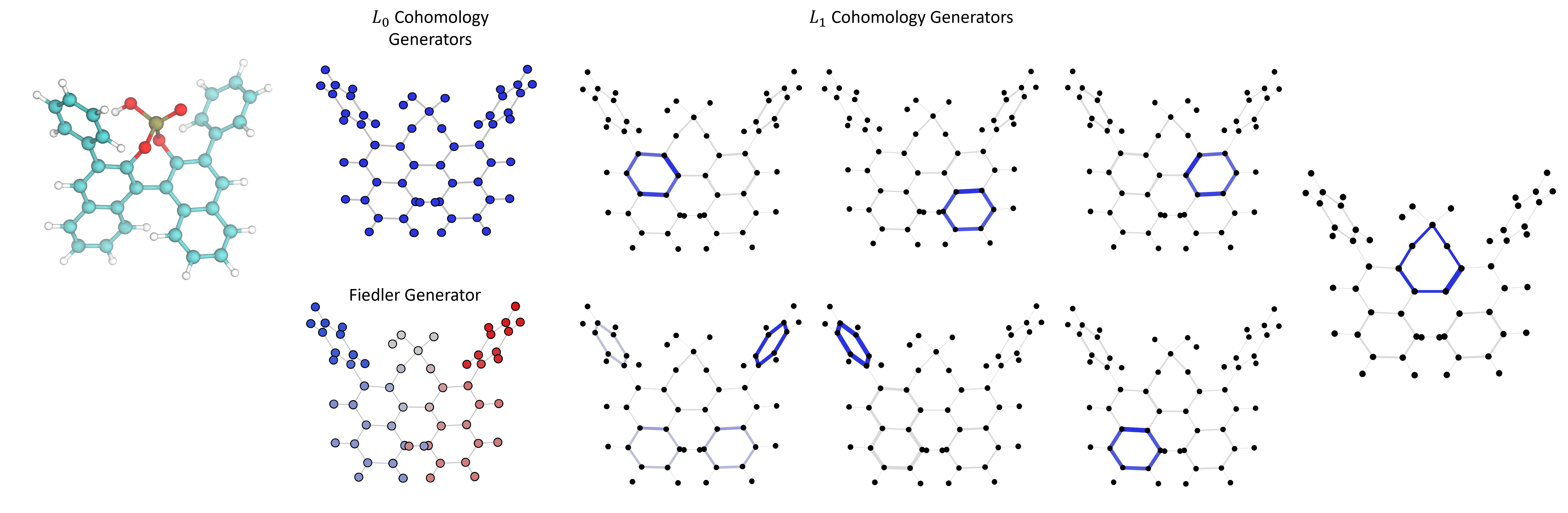}	\caption{Illustration of a catalyst 1\_{iv} with its $\mathbf{L}_0$ and $\mathbf{L}_1$ cohomology generators. The $\mathbf{L}_0$ cohomology generator shows all the $0$-simplices having equal value as it corresponds to the $\beta_0$, which represents the connected components of the catalyst. The $\mathbf{L}_1$ matrix has 7 cohomology generators in total, each representing a unique cocycle within the catalyst.}
	\label{fig:cycles}
\end{figure}

\begin{defn}[Cocycle distance]\label{defn:pdist2}
	Let $v$ and $w$ be two cohomology generators in $H^p(K)$ where $K$ is a simplicial complex. Their cocycle distance $d_s(v, w)$ is defined as the absolute difference in their $L_1$ norms: 
        \[
	d_s(v, w) = |\|v\|_{1}-\|w\|_{1}|.
	\]
\end{defn}
For illustration, consider the BINOL-phosphoramide \textbf{1\_iv}, a BINOL (1,1'-bi-2-naphthol)-based phosphoric acid catalyst. We construct the atomic level simplicial complex representation of the molecule, where each vertex represents an atom, and generate the corresponding Hodge Laplacian matrices $\mathbf{L}_0$ and $\mathbf{L}_1$. The cohomology generators from $\mathbf{L}_0$ and $\mathbf{L}_1$ are depicted in Figure \ref{fig:cycles}, where the color of each simplex $\sigma_i^p$ is darker in blue if its associated $|v_i|$ in the cohomology generator is larger. From $\mathbf{L}_0$, there is only one cohomology generator, representing a single connected component of the catalyst structure. Therefore, the cohomology generator in $\mathbf{L}_0$ assigns equal values to all vertices. On the other hand, $\mathbf{L}_1$ produces 7 cohomology generators in dark blue color in Figure \ref{fig:cycles} corresponding to unique cocycles. In addition, the non-cohomology generators associated with non-zero eigenvalues reveal local and global clustering patterns. The Fiedler vector effectively clusters the vertices into two groups, colored in blue and red. For example, Figure  \ref{fig:cycles} shows the Fiedler vector corresponding to the smallest non-zero eigenvalue of the Hodge Laplacian matrix, which clusters the vertices into two groups colored in blue and red.

\subsubsection*{Wasserstein Distance}
Recall from the beginning of this section that any cohomology generator $v$ in $H^p(K)$ is assumed to be normalized. Hence, given any two cohomology generators $v$ and $w$ in $H^p(K)$, we can square all the entries of $v$ and $w$ to obtain vectors $v'$ and $w'$, whose entries sum to 1,
\begin{equation}\label{eqn:1}
\text{i.e. }v' = \begin{blockarray}{cc}
\begin{block}{(c)c}
v_1^2 &  \sigma_{1}^p  \\
v_2^2 &  \sigma_{2}^p \\
\vdots &  \vdots \\
v_{n_p-1}^2 & \sigma_{n_p-1}^p  \\
v_{n_p}^2 &  \sigma_{n_p}^p \\
\end{block}
\end{blockarray}
\text{ and }\quad
w' = \begin{blockarray}{cc}
\begin{block}{(c)c}
w_1^2 &  \sigma_{1}^p  \\
w_2^2 &  \sigma_{2}^p \\
\vdots &  \vdots \\
w_{n_p-1}^2 & \sigma_{n_p-1}^p  \\
w_{n_p}^2 &  \sigma_{n_p}^p \\
\end{block}
\end{blockarray}.
\end{equation}
Since the vectors $v'$ and $w'$ have entries that sums to 1, this allows us to treat the values in $v'$ and $w'$ as two probability measures $m_1$ and $m_2$ respectively. Now, we can define the pairwise distance between any two cohomology generators $v$ and $w$ to be the Wasserstein distance \cite{kolouri2017optimal, villani2021topics} between $m_1$ and $m_2$.

\begin{defn}[Wasserstein distance]\label{defn:pdist3}
The probability measures $m_1$ and $m_2$ are the probability distributions obtained from the entries of $v'$ and $w'$ as defined in \eqref{eqn:1}. Let $\sigma_i^p$ and $\sigma_j^p$ be $p$-simplices in $K$. Then, $\xi(\sigma_i^p, \sigma_j^p)$ represents the amount of mass traveling from $\sigma_i^p$ to $\sigma_j^p$. We require that the transportation from $m_1$ to $m_2$ is mass-preserving, i.e., $\sum_{\sigma_j^p \in K} \xi(\sigma_i^p, \sigma_j^p) = m_1(\sigma_i^p)$ and $\sum_{\sigma_i^p\in K} \xi(\sigma_i^p, \sigma_j^p) = m_2(\sigma_j^p)$. The Wasserstein distance between $m_1$ and $m_2$, denoted by $W_1(m_1, m_2)$, is the minimum traveling distance that can be achieved, given by
	\begin{equation}
	W_1(m_1 , m_2) = \inf_\xi \sum_{\sigma_i^p\in K} \sum_{\sigma_j^p \in K} d(\sigma_i^p,\sigma_j^p)\xi(\sigma_i^p, \sigma_j^p).
	\end{equation}
The distance between two simplices, denoted by \(d(\sigma_i^p, \sigma_j^p)\), is defined as the minimum \(\ell^2\) distance between their corresponding vertices. Representing the sets of vertices in the simplices as $\sigma_i^p=\{x_1, x_2, \ldots, x_{p}\}$ and $\sigma_j^p=\{u_1, u_2, \ldots, u_{p}\}$, the distance is then given by
	\begin{equation}
	d(\sigma_i^p, \sigma_j^p) = \min_{\substack{1\leq i \leq p\\1\leq j\leq p}}||\mathbf{x}_i-\mathbf{u}_j||_2,
	\end{equation}
 where $\mathbf{x}_i$ and $\mathbf{u}_j$ refer to the coordinates of the vertices $x_i$ and $u_j$ in $\mathbb{R}^3$.
\end{defn}

\subsection*{Cohomology-based Gromov–Hausdorff Ultrametric Approach for Structural Similarity Measurement}
Given two molecular structures, $A$ and $B$, we construct simplicial complexes, denoted by $K_A$ and $K_B$ respectively, where vertices represent atoms. There are various approaches for constructing the simplicial complex \cite{carlsson2006algebraic}. For example, we can build a Vietoris-Rips complex such that a set of atoms is a simplex in $K$ if every atom pair in the set has a $\ell^2$ distance smaller than a specified filtration threshold. An alternative approach is to generate an Alpha complex, where a set of vertices forms a simplex when its filtration value is smaller than the threshold. For example, the filtration value of an atom pair is one-half of their $\ell^2$ distance, and the filtration value for a triplet of atoms is the radius of the circle that passes through it.

After constructing the simplicial complex representation of a molecular structure, cohomology generators can be computed from its 1-dimensional Hodge Laplacian matrix $\mathbf{L}_1$. Subsequently, we can compute a distance matrix comprised of pairwise distances between all cohomology generators from $\mathbf{L}_1$. Three types of distance matrices can be derived for each molecular structure, corresponding to Definition \ref{defn:pdist1}, \ref{defn:pdist2}, or \ref{defn:pdist3}.

In theory, we can compare the Gromov-Hausdorff distances between the metric spaces formed by the cohomology generators using a selected distance measure. However, it has been proven that computing \( d_{\text{GH}} \) is NP-hard \cite{schmiedl2017computational}. Therefore, we will adopt a more tractable ultrametric variation of the Gromov-Hausdorff distance, called the Gromov-Hausdorff ultrametric, which can be computed in polynomial time \cite{memoli2021gromov}. 

An ultrametric space $(X, d_X)$ is a metric space that satisfies the strong triangle inequality:
\[
\forall x, x', x''\in X, \text{ one has } d_X(x,x')\leq \max(d_X(x,x''),d_X(x'', x')).
\]
The Gromov-Hausdorff ultrametric \cite{zarichnyi2005gromov}, denoted by $u_{\text{GH}}$, measures distances between compact ultrametric spaces $X$ and $Y$ as follows.
\begin{defn}[The Gromov-Hausdorff ultrametric\cite{zarichnyi2005gromov}\label{defn:ugh}]
The Gromov-Hausdorff ultrametric $u_{\text{GH}}$ between \textit{compact} ultrametric spaces $X$ and $Y$ is
	\[
	u_{\text{GH}}(X,Y) := \inf d_H^Z(\varphi_X(X), \varphi_Y(Y)),
	\]
 where the infimum is taken over all \textit{ultrametric} spaces $Z$ and isometric embeddings $\varphi_X:X\hookrightarrow Z$ and $\varphi_Y:Y\hookrightarrow Z$.
\end{defn}

\begin{figure}[ht]
	\centering
	\includegraphics[width=\textwidth]{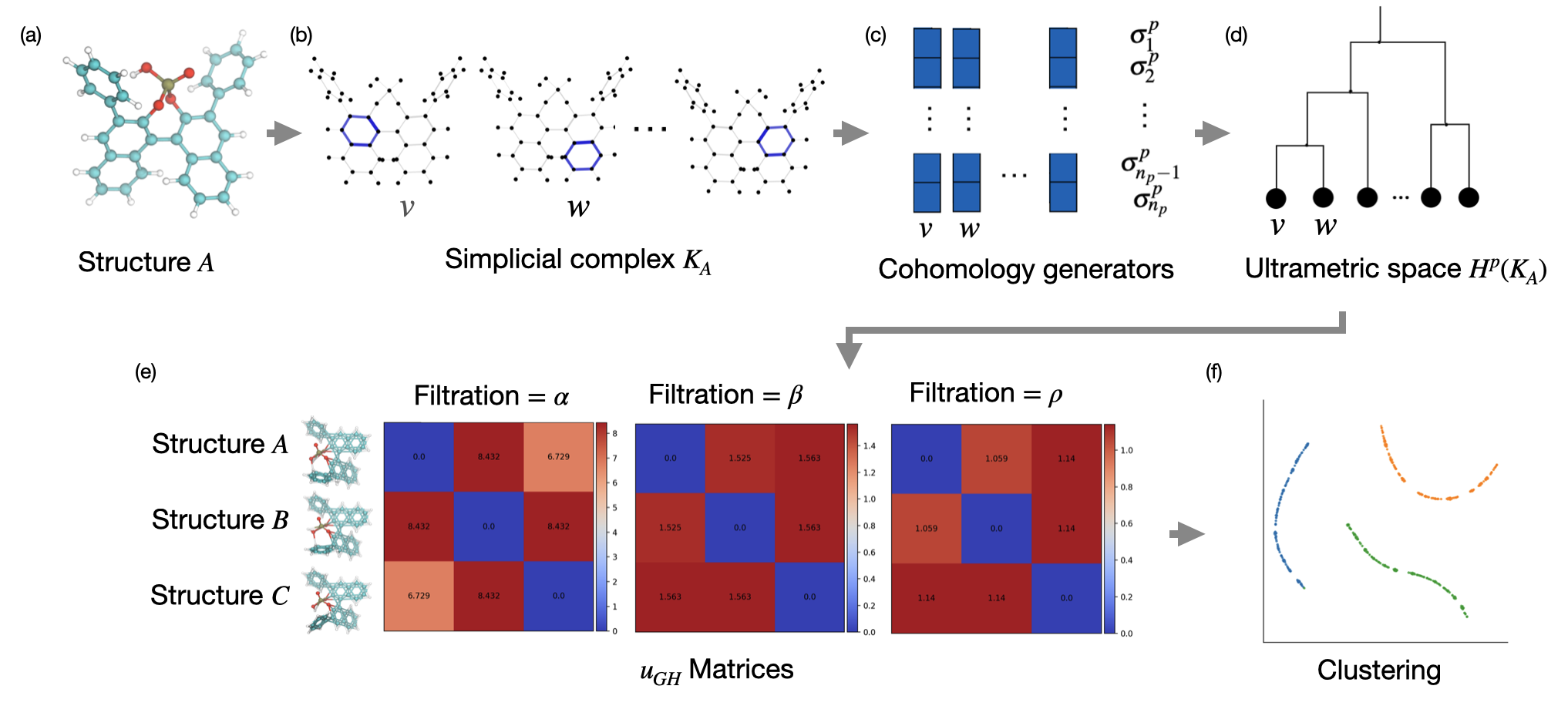}
	\caption{Schematic chart of the cohomology-based Gromov-Hausdorff approach for quantifying structural similarity. (a) Input coordinate data of structure $A$. (b) Construct the associated simplicial complex $K_A$. (c) Compute the cohomology generators. (d) Construct the ultrametric space $H^p(K_A)$ from the cohomology generators represented as a dendrograms. Nodes at the bottom correspond to cohomology generators and they merge at heights equal to their distances. (e) Compute the Gromov-Hausdorff ultrametric $u_{\text{GH}}$ between each pair of structures. Varying filtration values during the construction of the simplicial complex will result in multiple $u_{\text{GH}}$ matrices. (f) Use the $u_{\text{GH}}$ matrices as input to clustering algorithms.}

	\label{fig:flowchart}
\end{figure}

In our methodology, we transform metric spaces formed by the cohomology generators into an ultrametric space denoted by $H^1(K_A)$ and $H^1(K_B)$ using Algorithm 1 in \cite{memoli2021gromov}. This transformation allows us to compute the Gromov-Hausdorff ultrametric $u_{\text{GH}}(H^1(K_A), H^1(K_B))$ between two ultrametric spaces associated with molecular structures $A$ and $B$, thereby quantifying their structural similarity or facilitating clustering tasks.  An illustrative figure for the aforementioned workflow is presented in Figure \ref{fig:flowchart}.

\section*{Results}

\subsection*{Characterization of OIHP Material Structures}
\begin{figure}
	\centering
	\includegraphics[width=.5\textwidth]{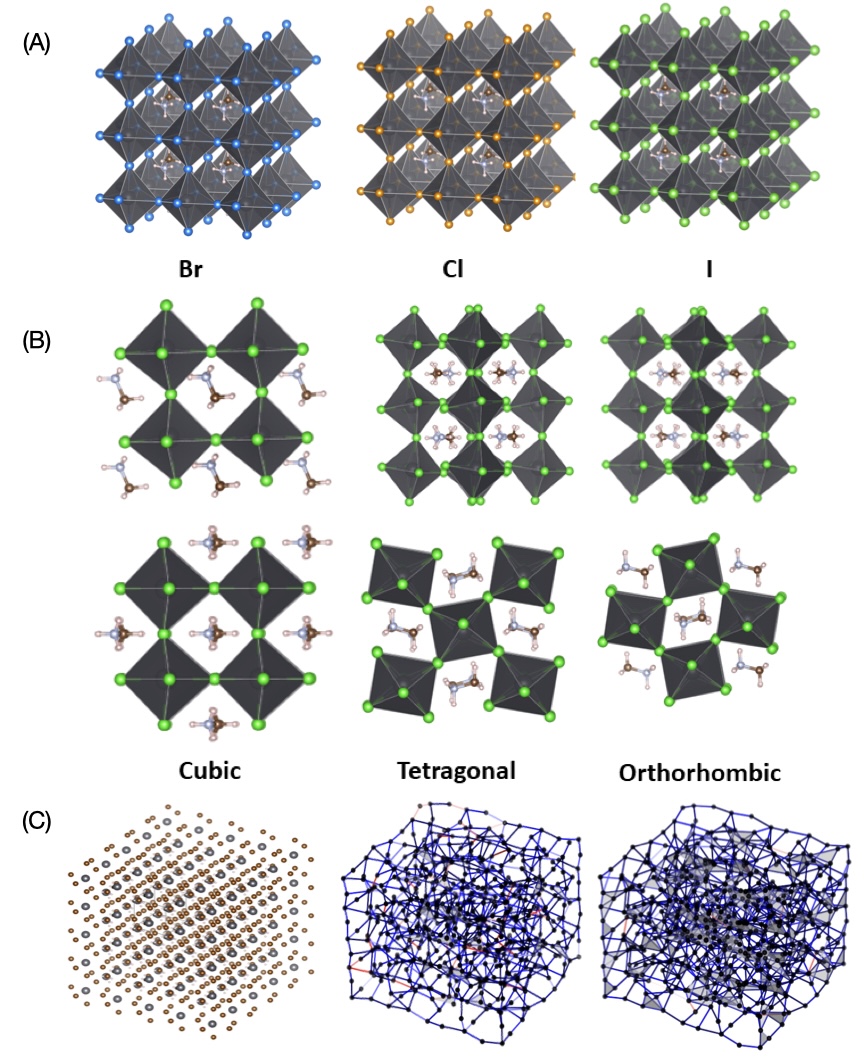}
	\caption{Illustration of 9 types of OIHP structures with formula MAPbX$_3$, where MA refers to Methylammonium, Pb is lead and X is bromine (Br), chlorine (Cl) or iodine (I). (a) Three OIHP structures with bromine, chlorine, and iodine. (b) Three phases of OIHP structures, i.e. cubic, orthorhombic, and tetragonal. In total, we consider all 9 possible combinations. (c) Illustration of the construction of the Alpha complex from MAPbBr$_3$ structures in the cubic phase (left) with a filtration value of 3.5~\AA ~(middle) and 4~\AA ~(right), where the edges are colored by the absolute value of a cohomology generator, with red corresponding to a larger absolute value. }
	\label{fig:perovskites}
\end{figure}

Organic-inorganic hybrid perovskite (OIHP) materials are favorable candidates for developing efficient and cost-effective solar cells. The stabilization of the OIHP structure is reliant on Van der Waals interactions and hydrogen bonding effects, which are closely associated with the distances between organic molecules and inorganic ions \cite{anand2022topological}. Our $u_{\text{GH}}$-based features can be employed for characterizing OIHP material structures. To demonstrate this, we examine the clustering of Methylammonium lead halides (MAPbX$_3$, X$=$Cl, Br, I). Three possible phases of MAPbX$_3$ are considered for each X-site atom, including the orthorhombic, tetragonal, and cubic phases as illustrated in Figure \ref{fig:perovskites} (a)-(b).  

We analyze 100 configurations from the molecular dynamics (MD) trajectories for each of the 9 OIHP structures, leading to 900 trajectories. Each MD trajectory arises from a molecular dynamics simulation to stabilize the initial configuration and consists of the 3D coordinates of atoms at a specific time. For each configuration, we construct an Alpha complex with 4 filtration thresholds of 3.5~\AA, 4~\AA, 5~\AA, and 6~\AA~  as shown in Figure \ref{fig:perovskites} (c) and calculate its 1-dimensional co-homology generators and associated pairwise distances using $L_1$ distance as defined in Definition \ref{defn:pdist1}. We chose to construct alpha complexes, which are based on the Delaunay triangulation of the points, because they better reflect the underlying geometry of atomic arrangement.

We consider the tasks of clustering 3 types of atoms for X-sites for each possible phase of MAPbX$_3$, by calculating the $u_{GH}$ between the 300 configurations of interest with the same phase for each of the 4 filtration values. This process generates a GH-based statistical feature vector with a length of $300 \times 4 = 1200$ for each configuration, allowing for differentiation among the three X-sites. We compare the K-means clustering by X-site atoms obtained using $u_{\text{GH}}$-based features with the ones achieved using the 3D coordinates, and two commonly used structural fingerprints: ECFP \cite{rogers2010extended} and MACC keys\cite{polton1982installation}. ECFP is a circular fingerprint that encodes atomic neighborhood patterns, while MACCS Keys captures the presence of specific chemical substructures. We chose to compare these features using only structural information without including additional information such as atomic number and weight, which could potentially bias the clustering algorithm towards identifying atom types. 

\begin{figure}[ht]
	\centering
	\includegraphics[width=\textwidth]{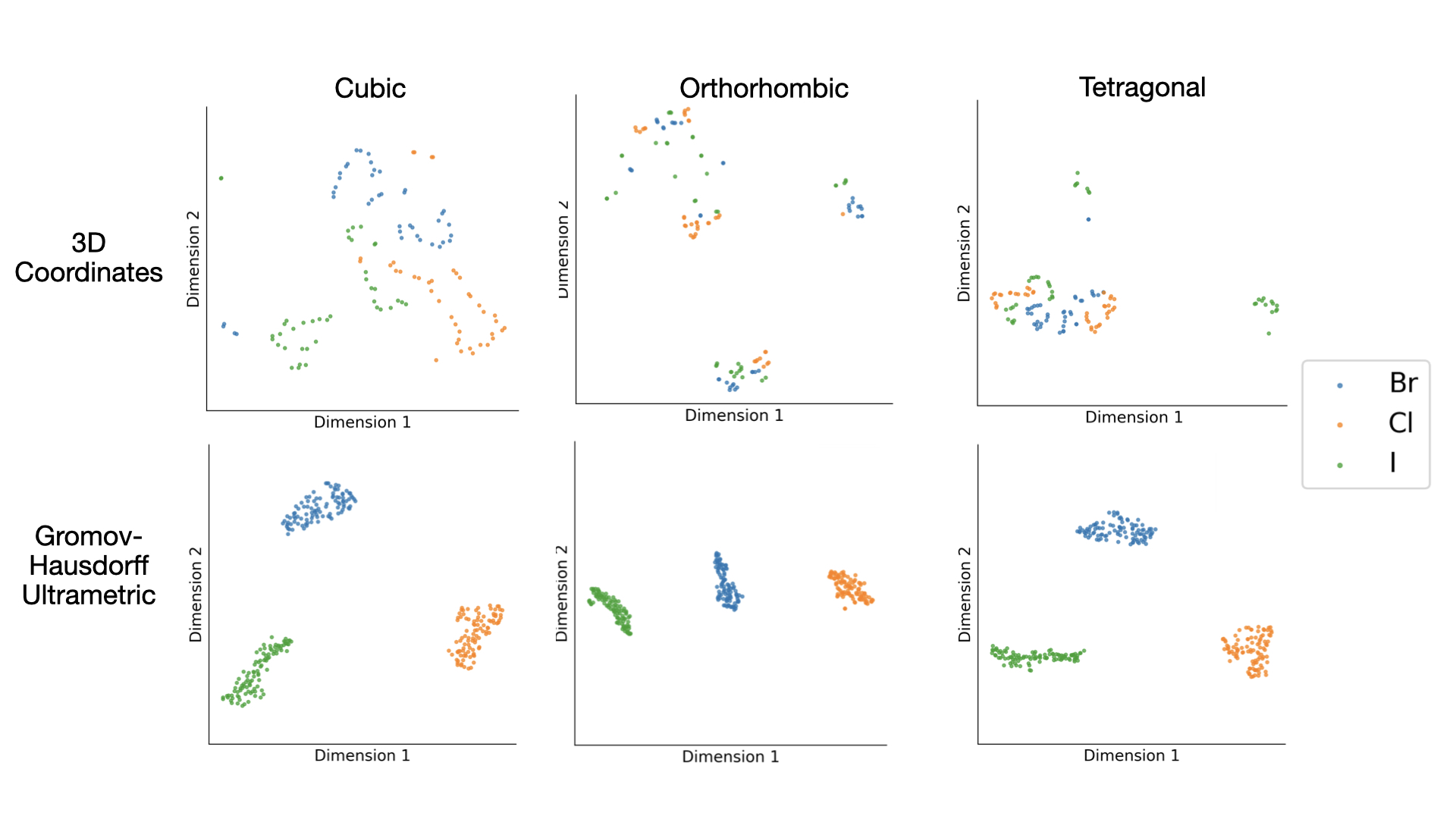}
	\caption{Visualization of OIHP molecular configurations using UMAP. The x-axis and y-axis represent the two dimensions after dimensionality reduction by the UMAP.}
	\label{fig:oihp_cluster}
\end{figure}

In Table \ref{table:ari}, we show the Adjusted Rand Index (ARI) of K-means clustering for the above methods, where 0 indicates a random assignment of clusters and 1 represents a perfect clustering. We also plotted low-dimensional visualization of these features using UMAP \cite{mcinnes2018umap} in Figure \ref{fig:oihp_cluster}, with the x-axis and y-axis representing the two dimensions after dimension reduction by UMAP. Directly using the 3D coordinates proves ineffective in distinguishing between various X-site atoms. Using the 3D coordinates alone proves ineffective in distinguishing between various X-site atoms as shown in low ARI. Incorporating some neighborhood information using ECFP improves the ARI but is still not ideal since the structures have similar neighborhoods pattern. MACC achieves decent clustering, although the ARI remains below 0.9 for the tetragonal structures. On the other hand, our Gromov-Hausdorff ultrametric-based methods exhibit  almost perfect clustering structures by X-site atoms. The UMAP plot in Figure \ref{fig:oihp_cluster} provides additional visual confirmation that our method achieves superior clustering in the low-dimensional projection of the data.

\begin{table}[h!]
\centering
\begin{tabular}{|c|c|c|c|c|}
\hline
Structure & 3D Coordinates & ECFP & MACC & Gromov-Hausdorff Ultrametric \\ \hline
Cubic  & 0.209  & 0.483  & 0.960 & \textbf{1.000}  \\ \hline
Orthorhombic & 0.009  & 0.030  & \textbf{1.000} & \textbf{1.000}  \\ \hline
Tetragonal  &  0.436  & 0.322  & 0.827 & \textbf{0.980}  \\ \hline
\end{tabular}
\caption{Adjusted Rand Index of clustering based on different features}
	\label{table:ari}

\end{table}

\section*{Discussion}

In this work, we propose a new framework that employs the cohomology-based Gromov-Hausdorff ultrametric ($u_{GH}$) for quantifying structural similarity between molecular structures. In our workflow, we build simplicial complex representations of molecular structures and compute the $u_{GH}$ between their respective cohomology generator spaces. The cohomology generators effectively encode the local topological invariants, revealing cyclic patterns formed by edges. We illustrate the application of the cohomology-based Gromov-Hausdorff distance using organic-inorganic halide perovskites (OIHP) data. In our numerical experiments, the $u_{GH}$-based approach demonstrated effectiveness in clustering various structures, achieving the highest Adjusted Rand Index compared to other structural features tested including the 3D coordinates, ECFP, and MACC.

We demonstrate the application of the cohomology-based $u_{GH}$ approach in quantifying structural similarities. However, there are some limitations to our current work. For instance, we focused solely on the cohomology generators of the first-order Hodge Laplacian, which capture loop patterns. In the future, exploring non-cohomology generators may reveal clustering information, and higher-order Hodge Laplacians could unveil more complex structures such as cavities. Another limitation is that the choice of kernel vectors for the Hodge Laplacian is not unique. Therefore, for the same structure, we may have different choices for cohomology generators. It remains an open question for future research to determine how to find the optimal cohomology generators that minimize the $u_{GH}$ between two atoms. Other potential future directions involve incorporating cohomology-based $u_{GH}$ features into machine learning models for structure design and prediction. Topological deep learning models have also demonstrated great potential in molecular property prediction and analysis \cite{ballester2024attending}. Additionally, our numerical experiments are limited to small molecules from chemistry and physics with hundreds of atoms. A promising avenue involves applying the proposed method to larger biological molecules with tens of thousands of molecules. For instance, it could be employed to quantify structural similarities between protein structures in drug design \cite{dean2012molecular}. 
\bibliography{main}% Produces the bibliography via BibTeX.
 
%\section*{Bibliography}

\section*{Acknowledgements}

This work was supported in part by the Singapore Ministry of Education Academic Research fund Tier 1 grant RG16/23, Tier 2 grants MOE-T2EP20120-0010 and MOE-T2EP20221-0003. X.G. acknowledges support from NTU Presidential Postdoctoral Fellowship 023545-00001. The authors would like to thank Yipeng Zhang for his contributions to the revision of this manuscript, and Chuan-Shen Hu for his helpful  suggestions.

\section*{Author contributions statement}

W.T. and K.X. designed the experiments. J.W. and X.G. conducted the experiments and analyzed the results. J.W. drafted the initial manuscript. X.G. revised the manuscript. All authors reviewed and contributed to the final manuscript.

\section*{Additional information}

The authors declare that there are no competing interests related to this work.

\section*{Data availability}
Code and data are available in GitHub repository https://github.com/XueGong-git/pyGH.

\end{document}